\documentclass[sn-mathphys-num]{sn-jnl}

\usepackage{graphicx}%
\usepackage{multirow}%
\usepackage{amsmath,amssymb,amsfonts}%
\usepackage{amsthm}%
\usepackage{mathrsfs}%
\usepackage[title]{appendix}%
\usepackage{xcolor}%
\usepackage{textcomp}%
\usepackage{manyfoot}%
\usepackage{booktabs}%
\usepackage{listings}%

\theoremstyle{thmstyleone}%
\newtheorem{theorem}{Theorem}
\theoremstyle{thmstyletwo}%
\newtheorem{lem}{Lemma}[section]
\theoremstyle{thmstylethree}%

\newcommand{\dE}{\mathbb{E}}

\newcommand{\cN}{\mathcal{N}}

\begin{document}

\title[Lazy reinforced random walk]{On the asymptotics of a lazy reinforced random walk}

\author[1]{\fnm{Manuel} \sur{Gonz\'alez-Navarrete}}\email{manuel.gonzaleznavarrete@ufrontera.cl}

\author[2]{\fnm{Rodrigo} \sur{Lambert}}\email{rodrigolambert@yahoo.com.br}
\equalcont{These authors contributed equally to this work.}

\author*[3]{\fnm{V\'ictor Hugo} \sur{V\'azquez Guevara}}\email{victor.vazquezg@correo.buap.mx}
\equalcont{These authors contributed equally to this work.}

\affil[1]{\orgdiv{Departamento de Matem\'atica y Estad\'istica}, \orgname{Universidad de La Frontera}, \orgaddress{\street{Avda. Francisco Salazar 1145}, \city{Temuco}, \postcode{4811230}, \state{La Araucan\'ia}, \country{Chile}}}

\affil[2]{\orgdiv{Faculdade de Matem\'atica}, \orgname{Universidade Federal de Uberl\^andia}, \orgaddress{\street{Av. Jo\~ao Naves de \'Avila, 2121 - Santa M\^onica}, \city{Uberl\^andia}, \postcode{38408100}, \state{Minas Gerais}, \country{Brasil}}}

\affil*[3]{\orgdiv{Facultad de Ciencias F\'isico Matem\'aticas}, \orgname{Benem\'erita Universidad Aut\'onoma de Puebla}, \orgaddress{\street{San Claudio y R\'io Verde	}, \city{Puebla}, \postcode{72570}, \state{Puebla}, \country{M\'exico}}}


\abstract{Based on a martingale theory approach, we present a complete characterization of the asymptotic behaviour of a lazy reinforced random walk (LRRW) which shows three different regimes (diffusive, critical and superdiffusive). This allows us to prove versions of the law of large numbers, the quadratic strong law, the law of iterated logarithm, the almost sure central limit theorem and the functional central limit theorem in the diffusive and critical regimes. In the superdiffusive regime we obtain a strong convergence to a random variable, including a central limit theorem and a law of iterated logarithm for the fluctuations.}




\maketitle

\section{Introduction}

\hspace{.5cm} Over the years, the study of stochastic processes with infinite memory has received special attention in literature. The idea that the whole past of the process has influence on its behaviour and asymptotics has been investigated by many researchers in mathematics, statistics and physics. 

A very well-known family of processes which may present infinite memory is the one of P\'olya urn models. Depending on the nature of the urn, its initial condition may either determine strong convergence (a trivial example is when the urn begins with only one color of balls and the dynamics doesn't allow to add new ones) or provide the parameters of such limiting distribution (for instance, the classical P\'olya-Eggenberger model converges weakly to a beta random variable with parameters that depend on the initial configuration of the urn). We refer the reader to \cite{Mahmoud} for a complete survey on P\'olya urns, and to \cite{Janson} for a comprehensive list of asymptotic results and examples for such model (and also its links with branching processes).
 
%

Going further on infinite memory models, a remarkable example appeared in 2004. It describes a random walk in which at each step the walker draws uniformly a time of the past and chooses randomly if he will repeat its own movement on such time or will do the opposite. Inspired by the common belief that elephants have very good memory, this process was baptised as the elephant random walk (ERW). It was proposed in \cite{Schutz} and its asymptotic characterization was completed in \cite{ERWBercu,Lucile1,Lucile2,kubota2019gaussian,gue}. Several variations of the ERW in which only a portion of the past (in many different ways) is taken into account were proposed in \cite{Gut}. In a different way, \cite{gonzalez2018non,gonzalez2019diffusion} introduced the possibility that at each step the walker may either base its current movement on its whole past or act independently. In the sequel, a reinforced random walk with the so-called memory lapse property was proposed in \cite{GH}. By studying its asymptotics (law of large numbers, central limit theorem, law of iterated logarithm), it was shown the dependency of the limit quantities and the memory lapse parameter (which guides to the limit analysis into the diffusive, critical and superdiffusive regimes).

In this paper we complement the asymptotical analysis of the reinforced random walk discussed in \cite{GH} by considering a more natural martingale which will lead us to obtain refinements of the earlier results. In addition, new findings of the almost sure central limit theorem, functional central limit theorem and convergence to even moments of Gaussian distribution type will be explored.



The rest of paper is organized as follows: the following section introduces the random walk model we will deal with and will also present the associated martingale which will become the key in order of achieving the proofs of all the discussed results. Section 3 states the main results, which deals with the asymptotics of the random walk under consideration. Last section contains the proofs of our results as well as an initial martingale analysis of the LRRW.

\section{The Random Walk}
\hspace{.5cm}This section is devoted to introduce the lazy reinforced random walk (LRRW), which is the heart of this work. It also presents some notation and tools that will be very useful along the paper.

The random walk we will deal with, whose movement and position at time $n\geq 1$ are denoted by $X_n$ and $S_n$ respectively, behaves as follow: 

\begin{enumerate}
\item At time $n=0$ the LRRW is at the origin; i.e., $S_0=0$.
\item At time $n=1$, we have that $S_1=X_1$ has the following distribution 
$$         X_1 = \left\{
	       \begin{matrix}
			1,& \text{with probability $p$,}\\
		-1,&\text{\text{with probability $q$,}}\\
		0,&\text{with probability $r$.}
		 
	  \end{matrix}%
	\right.
  $$
\item For each $n\geq1$  set
\begin{equation}\label{RWdef}
X_{n+1}=Y_n \alpha_{n+1}X_{U_n}+\left(1-Y_n\right)\alpha_{n+1},
\end{equation}%
where $U_n$ is a discrete uniform random variable on $\left\{1,2,\ldots,n\right\}$, $Y_n$ posses the Bernoulli distribution with parameter $\theta\in[0,1)$ and
$$         \alpha_{n+1} = \left\{
	       \begin{matrix}
			1,& \text{with probability $p$,}\\
		-1,&\text{\text{with probability $q$,}}\\
		0,&\text{with probability $r$.}
		 
	  \end{matrix}%
	\right.
  $$
with $p+q+r=1$. In addition, we assume that $\alpha_n$, and $U_n$ are independent and that $Y_n$ is independent of the LRRW's past. Then, the position of the LRRW at time $n\geq 0$ is given by
\begin{equation}\label{Sn}
S_{n+1}=S_n+X_{n+1}.
\end{equation}

\end{enumerate}

At this point, since $X_{n+1}$ depends on $Y_n$, $\alpha_{n+1}$ and $X_{U_n}$, lets take a break to discuss the rule played by each term in \eqref{RWdef}. 

First, if $Y_n=0$, we can see that the step will not depend on its past. Nevertheless, at the next time that $Y_k=1$  (for some $k>n$), the past will be taken into account (unless $\alpha_{k+1}=0$) by the walker. In other words, the random variable $Y_n$ turns on (or off) the memory of the walker. We call this the \emph{memory lapse property}. Moreover, the random variables $\alpha_{n+1}$ say essentially if the walker won't walk ($\alpha_{n+1}=0$), repeat ($\alpha_{n+1}=1$), or do the opposite ($\alpha_{n+1}=-1$) that she (he) has done at step $U_n$. Therefore, the $\alpha_n$ process plays an agreement (or disagreement) rule.

Finally, note that, for $\theta=1$, we obtain the already characterized model in \cite{GS2}.

Now we will derive some tools (and present some additional notation).
If $\left(\mathcal{F}_n\right)$ is the increasing sequence of $\sigma$-algebras $\mathcal{F}_n=\sigma\left(X_1,\ldots,X_n\right)$; i.e. $\mathcal{F}_n$ represents the knowledge on the process $(X_i)$ up to time $n$, we get 
\begin{equation}\label{e1}
\mathbb{E}\left[X_{n+1}|\mathcal{F}_n\right]=\mathbb{E}[Y_n]\mathbb{E}[\alpha_{n+1}]\mathbb{E}[X_{U_n}|\mathcal{F}_n]+\mathbb{E}[\alpha_{n+1}] \left(1-\mathbb{E}[Y_n] \right) \hspace{.5cm}\text{a.s.,}
\end{equation}
By using the definition of $S_n$ combined with the total probability theorem we find that 
$$\mathbb{E}[X_{U_n}|\mathcal{F}_n]=\frac{S_n}{n}\hspace{.5cm}\text{a.s.}$$
Combining the above expression with \eqref{e1} we conclude that
\begin{equation}\label{expmov}
\mathbb{E}\left[X_{n+1}|\mathcal{F}_n\right]=\alpha\frac{S_n}{n}+\omega \hspace{.5cm}\text{a.s.}
\end{equation}
where $\alpha:=\theta(p-q)$ and $\omega:=(p-q)(1-\theta)$. Therefore
\begin{equation}\label{expos}
\mathbb{E}\left[S_{n+1}|\mathcal{F}_n\right]=\gamma_n S_n+\omega \hspace{.5cm}\text{a.s.,}
\end{equation}%
where $\gamma_n=1+\frac{\alpha}{n}.$
Now it is time to present our main ingredient for the proofs. It is the associated martingale, that will lead us on the asymptotic analysis of the LRRW. Let the sequence $(M_n)$, given by $M_0=0$ and for $n\geq 1$  by
\begin{equation} \label{martingale}
M_n=a_n S_n-\omega A_n.
\end{equation}
In what follows, we will define $a_n$ and $A_n$. First of all, if $\Gamma$ stands for the Euler gamma function, then
\begin{equation}\label{an}
a_n=\prod_{k=1}^{n-1}\gamma_k^{-1}=\frac{\Gamma(n)\Gamma(\alpha+1)}{\Gamma(n+\alpha)}\sim \frac{\Gamma(1+\alpha) }{n^\alpha}.
\end{equation}%
We also recall that the equality above has been obtained by using recursively the property $\Gamma(\alpha+1)=\alpha\Gamma(\alpha)$, and the approximation can be obtained by using Stirling's approximation for the Euler gamma function. Now we define $A_0=0$, and for $n\geq 1$
\begin{equation}\label{bigAn}
A_n=\sum_{k=1}^n a_k.
\end{equation}
Additionally, we observe from \eqref{expos} and \eqref{an} that for almost every realization it holds
\begin{eqnarray*}
\mathbb{E}[M_{n+1}|\mathcal{F}_n]&=&a_{n+1}(\gamma_n S_n+\omega)-\omega A_{n+1}\\
&=&a_n S_n-\omega A_n \\ 
&=& M_n \ .
\end{eqnarray*}
Which proves that $(M_n)$ is a discrete time martingale with respect to the filtration $\left(\mathcal{F}_n\right)$.

Asymptotic behaviour of the LRRW posses different attributes according with the following nomenclature:
\begin{enumerate}
\item If $\alpha<1/2$ it is said that the LRRW is in the diffusive regime.
\item If $\alpha=1/2$ then the regime is called critical, and
\item If $\alpha>1/2$ then the regime is labeled as superdiffusive.
\end{enumerate}
Hence, we will conduct the analysis of the LRRW according to the previous classification.


\section{Main Results}

In all the sequel, we will use the following notation in the asymptotic analysis of the LRRW:
\begin{equation}
\label{notation}
\tau=(1-\theta)(p+q), \ \gamma=\theta(p+q) \text{ and } \sigma^2 = \frac{\tau}{1-\gamma}-\left(\frac{\omega}{1-\alpha}\right)^2
\end{equation}
where $\alpha$ and $\omega$ as in \eqref{expmov}.

\subsection{The diffusive regime}

We initiate by studying the asymptotics of the LRRW in the diffusive regime. The following result deals with the law of large numbers and its corresponding convergence rate. 

\begin{theorem}\label{CVG1}
If $\alpha<1/2$, then we have the following almost sure convergence 
\begin{equation}\label{as1}
\lim_{n \rightarrow \infty} \frac{S_n}{n}=\frac{\omega}{1-\alpha},
\end{equation}%
to be precise
\begin{equation} \label{as00}
\left( \frac{S_n}{n}-\frac{\omega}{1-\alpha} \right)^2 =O\left( \frac{\log n}{n}\right), \hspace{.2cm}\text{a.s.}
\end{equation}
\end{theorem}

In addition, we have the convergence to the even moments of Gaussian distribution.

\begin{theorem}
\label{thmmoments}
If $\alpha < 1/2$ then the following almost sure convergence holds, as $n \to \infty$
\begin{equation}
\label{moments1}
\begin{array}{lll}
  \displaystyle\frac{1}{\log n}\sum_{k=1}^n k^{r-1}\left(\frac{S_k}{k} - \frac{\omega}{1-\alpha} \right)^{2r} \rightarrow \frac{(\sigma^2)^r(2r)!}{2^rr!(1-2\alpha)^r }
\end{array}
\end{equation}
\end{theorem}
In particular, by replacing $r=1$ in \eqref{moments1}, we obtain the quadratic strong law.

We also have the law of iterated logarithm 
\begin{theorem}\label{CVG3}
If $\alpha<1/2$, then 
\begin{eqnarray*}
&\displaystyle\limsup_{n \rightarrow \infty}&  \left(\frac{n}{2 \log \log n}\right)^{1/2} \left(\frac{S_n}{n}-\frac{\omega}{1-\alpha}\right) \\
&=& -\liminf_{n \rightarrow \infty}  \left(\frac{n}{2 \log \log n}\right)^{1/2} \left(\frac{S_n}{n}-\frac{\omega}{1-\alpha}\right) \\
&=& \frac{\sigma}{\sqrt{1-2\alpha}} \hspace{.5cm}\text{a.s.}
\end{eqnarray*}
\end{theorem}

One version of the almost sure central limit theorem is now enunciated for the diffusive regime:
\begin{theorem}
\label{thmASCLT}
If $\alpha <1/2$ then we have the following almost sure convergence of empirical measures, as $n \to \infty$
\begin{equation} \label{asclt}
  \displaystyle\frac{1}{\log n}\sum_{k=1}^n \frac{1}{k} \delta_{\sqrt{k}\left(\frac{S_k}{k} - \frac{\omega }{1-\alpha}\right)} \Rightarrow G  \ \  \text{ a.s}
\end{equation}
where $\delta_{x}(A)$ equals 1, if $x \in A$, or 0 otherwise. $G$ is the Gaussian measure $N(0, \sigma^2 /(1-2 \alpha))$ and $\Longrightarrow$ stands for convergence in distribution. 
\end{theorem}

Finally, we focus our attention on the functional convergence of the LRRW. Then, denote by $D([0,\infty[)$ the Skorokhod space of right-continuous functions with left-hand limits (c\`adl\`ag functions).

\begin{theorem}
\label{T-FCLT-DR}
If $\alpha <1/2$, we have the distributional convergence in $D([0,\infty[)$, as $n \to \infty$
\begin{equation}
\label{FCLT-DR}
\left( \sqrt{n}\Big(\frac{S_{\lfloor nt \rfloor}}{\lfloor nt \rfloor}-\frac{\omega}{1-\alpha}\Big), t \geq 0\right) \Longrightarrow \big( W_t, t \geq 0 \big)
\end{equation}
where $\big( W_t, t \geq 0 \big)$ is a real-valued centred Gaussian process starting at the origin with covariance given, for all $0<s \leq t$, by
$$
\dE[W_s W_t]= \frac{\sigma^2}{(1-2\alpha)t} \Big(\frac{t}{s}\Bigr)^{\alpha}.
$$
In particular, if $t=1$, we have the asymptotic normality (which has already been demonstrated in Theorem 2 of \cite{GH}).
\end{theorem}

\subsection{The critical regime}

In the present section, we will investigate the corresponding asymptotic results for the LRRW in the critical regime. We show the law of large numbers and its convergence rate.

\begin{theorem}\label{CVG11}
If $\alpha = 1/2$, then we have the following almost sure convergence 
\begin{equation}\label{as11}
\lim_{n \rightarrow \infty} \frac{S_n}{n}=2\omega,
\end{equation}
more precisely
\begin{equation} \label{as001}
\left( \frac{S_n}{n}-2\omega \right)^2 =O\left( \frac{\log n \log \log n}{n}\right), \hspace{.2cm}\text{a.s.}
\end{equation}
\end{theorem}

In addition, the corresponding convergence to even moments of Gaussian distribution is as follows:

\begin{theorem}
\label{thm4}
If $\alpha = 1/2$, then the following almost sure convergence holds, as $n \to \infty$
\begin{equation}
\label{moments2}
\begin{array}{lll}
  \displaystyle\frac{1}{\log \log n}\sum_{k=1}^n \left(\frac{1}{k \log k}\right)^{r+1} k^{r-1} \left(\frac{S_k}{k} - 2\omega\right)^{2r}\rightarrow \frac{(\sigma^2)^r(2r)!}{2^r r!}
\end{array}
\end{equation}
\end{theorem}
Which implies that, in the particular case $r=1$, the quadratic strong law holds.
%
%

The law of iterated logarithm is as follows
\begin{theorem}\label{CVG31}
If $\alpha = 1/2$, then the following almost sure convergence holds
\begin{eqnarray*}
\limsup_{n \rightarrow \infty}  \left(\frac{n}{2 \log n \log \log \log n}\right)^{1/2} \left(\frac{S_n}{n}-2\omega\right) \\
= -\liminf_{n \rightarrow \infty}  \left(\frac{n}{2 \log n \log \log \log n}\right)^{1/2} \left(\frac{S_n}{n}-2\omega\right)\\
 = \sqrt{\frac{\tau}{1-\gamma}-4\omega^2} \hspace{.5cm}\text{a.s.}
\end{eqnarray*}
\end{theorem}

\begin{theorem}
\label{thmasclt}
If $\alpha = 1/2$ we have the almost sure convergence, as $n \to \infty$
\begin{equation} \label{asclt2}
  \displaystyle\frac{1}{\log \log n}\sum_{k=1}^n \frac{1}{k \log k} \delta_{\sqrt{\frac{k}{\log k}}\left(\frac{S_k}{k} - 2\omega \right)} \Rightarrow G  \ \  \text{ a.s.,}
\end{equation}
where $G$ stands for the Gaussian measure $N(0,\sigma^2)$.

\end{theorem}

Our next result deals with the proper version of the functional central limit theorem in the critical regime.

\begin{theorem}
\label{T-FCLT-CR}
If $\alpha = 1/2$, we have the distributional convergence in $D([0,\infty[)$, as $n \to \infty$
\begin{equation}
\label{FCLT-CR}
\left( \sqrt{\frac{n^t}{\log n}}\Big(\frac{S_{\lfloor n^t \rfloor}}{\lfloor n^t \rfloor}-2\omega\Big), t \geq 0\right) \Longrightarrow \left( \left( \frac{\tau}{1-\gamma}-4\omega^2 \right) B_t, t \geq 0 \right)
\end{equation}
where $\big( B_t, t \geq 0 \big)$ is a standard Brownian motion.
In particular, we have the asymptotic normality for $t=1$ (already proved in Theorem 2 of \cite{GH}).
\end{theorem}

\subsection{The superdiffusive regime}

The superdiffusive regime leads us to a more exhaustive analysis which is, in some way, apart from the diffusive and critical regimes.

\begin{theorem}
\label{T-ASP-SR}
If $\alpha> 1/2$, we have the almost sure convergence, as $n \to \infty$
\begin{equation}
\label{FSLLN-SR}
\left( n^{1-\alpha}\Big(\frac{S_{\lfloor nt \rfloor}}{\lfloor nt \rfloor}-\frac{\omega}{1-\alpha}\Big), t > 0\right) \longrightarrow \Big( \frac{1}{t^{1-a}}L, t > 0 \Big)
\end{equation}
where $L$ is a non-degenerated random variable such that
\begin{equation}\label{meanL}
\mathbb{E}[L]=\frac{\beta(1-\alpha)-\omega}{\Gamma(\alpha+1)(1-\alpha)}
\end{equation}
where $\beta:= p-q$, and
\begin{equation} \label{varL}
\mathbb{E}[L^2]=\frac{\nabla}{\Gamma(2\alpha+1)}+2\omega\left(\frac{1}{(1-\alpha)\Gamma(\alpha)}\right)^2
\end{equation}
where $\nabla:=p+q+\frac{\tau}{(1-\gamma)(2\alpha-1)}-\frac{2\alpha \omega^2}{(2\alpha-1)(\alpha-1)^2} +4\left[\frac{\omega \alpha(\beta-1)}{(\alpha-1)^2}\right]+\frac{r\gamma^2}{2\alpha-\gamma}.$
\end{theorem}

%
%

Finally, we analyze the limit distribution of the LRRW

\begin{theorem}
\label{Gaussfluct}
If $\alpha>1/2$, then
\begin{equation}
    \label{T_cnv_w}
    \sqrt{n^{2\alpha-1}}\left( n^{1-\alpha}\Big(\frac{S_n}{n}-\frac{\omega}{1-\alpha}\Big)-L \right) 
\Longrightarrow \cN\Big(0,\frac{\sigma^2}{2\alpha-1}\Big)
 \text{ as }  \ n \to \infty
\end{equation}
   and
    \begin{equation}
    \label{LIL3}
        \limsup_{n \to \infty} \pm \frac{\sqrt{n^{2\alpha-1}}\left( n^{1-\alpha}\Big(\frac{S_n}{n}-\frac{\omega}{1-\alpha}\Big)-L \right)}{\sqrt{ \log \log n}} = \sqrt{\frac{2\sigma^2}{2\alpha-1}} \ \text{ a.s }
    \end{equation}
\end{theorem}

\section{Demonstrations of achieved results}

\subsection{Martingale Analysis}\label{secA1}

From the definition of the proposed martingale given in \eqref{martingale}, we observe that

\begin{eqnarray}
\Delta M_{n}&=&M_{n}-M_{n-1}=a_{n}\left(S_{n}-\gamma_{n-1} S_{n-1}-\omega\right)\\
&=&a_{n} \xi_{n}, 
\end{eqnarray}
where; for $n\geq 1$, $\xi_n:=S_n-\mathbb{E}[S_n|\mathcal{F}_{n-1}]=S_{n}-(\omega+\gamma_{n-1} S_{n-1})$ then
\begin{equation} \label{martingale2}
M_n=\sum_{k=1}^n a_k \xi_k
\end{equation}

In addition, equation \eqref{expos}, leads us to see that
\begin{equation} \label{epsn}
\mathbb{E}\left[\xi_{n+1}|\mathcal{F}_n\right]=0\hspace{.2cm}\text{a.s.,}\\
\end{equation}

From the total probability theorem, we have that, almost surely
$$\mathbb{E}[X_{U_n}^2|\mathcal{F}_{n}]=\frac{Z_n}{n},$$
where 
\begin{equation}\label{Zn}
Z_n=\sum_{k=1}^n 1_{\left\{1\right\}}\left(X_k\right)+\sum_{k=1}^n 1_{\left\{-1\right\}}\left(X_k\right)=\sum_{k=1}^n X_k^2
\end{equation}

that together lead us to
\begin{eqnarray}
\mathbb{E}[X_{n+1}^2|\mathcal{F}_{n}]&=&\mathbb{E}[ \left(Y_n \alpha_{n+1}X_{U_n}+\left(1-Y_n\right)\alpha_{n+1} \right)^2 |\mathcal{F}_{n}] \notag \\
&=&\gamma\frac{Z_n}{n}+\tau \hspace{.2cm}\text{a.s.,} \label{xsqq}
\end{eqnarray}
where $\gamma=(p+q)\theta$ and $\tau=(1-\theta)(p+q)$. From this, we find that, a.s.
\begin{equation}\label{Snsqq}
\mathbb{E}[S_{n+1}^2|\mathcal{F}_{n}]=S_n^2+2S_n\left(\frac{\alpha S_n}{n}+\omega\right)+\gamma\frac{Z_n}{n}+\tau.
\end{equation}
%

Therefore, from \eqref{epsn} we have that

\begin{equation} \label{segunda}
\mathbb{E}\left[\xi^2_{n+1}|\mathcal{F}_n\right]=\gamma\frac{Z_n}{n}+\tau-\left(\frac{\alpha S_n}{n}+\omega \right)^2,
\end{equation}
and given that since $Z_n\leq n$, we find that
\begin{equation} \label{cota112}
\mathbb{E}\left[\xi^2_{n+1}|\mathcal{F}_n\right]  \leq  \gamma+\tau.
\end{equation}
 This fact implies that
\begin{equation}\label{cota2}
\sup_{n\geq 0} \mathbb{E}\left[\xi^2_{n+1}|\mathcal{F}_n\right] < \infty.
\end{equation}

On the same direction it can be shown that

\begin{equation}
\mathbb{E}[S_{n+1}^3|\mathcal{F}_n]=S_n^3+\left[\alpha\frac{S_n}{n}+\omega\right]\left[3S_n^2+1\right]+3S_n\left[\gamma\frac{Z_n}{n}+\tau\right] \ \ \text{a.s.}
\end{equation}
And also

\begin{equation}
\mathbb{E}[S_{n+1}^4|\mathcal{F}_n]=S_n^4+\left(\omega+\alpha\frac{S_n}{n}\right)\left[4S_n^3+4S_n\right]+\left[\gamma \frac{Z_n}{n}+\tau \right]\left[6S_n^2+1\right] \hspace{.2cm}\text{a.s.}
\end{equation}
and hence, from the definition of $a_n$, we conclude that almost surely

\begin{eqnarray}
\mathbb{E}[\xi_{n+1}^4|\mathcal{F}_n]&=&\left(\tau+\gamma\frac{Z_n}{n}\right)+\left(\omega+\alpha\frac{S_n}{n}\right)^2\left[6\left(\tau+\gamma\frac{Z_n}{n}\right)-4\right]\\
&-&3\left(\omega+\alpha\frac{S_n}{n}\right)^4,\notag
\end{eqnarray}
nevertheless, as $|S_n|$ and $Z_n$ are both lower or equal than $n$, we conclude that

\begin{equation}\label{cuarta}
\sup_{n\geq 0} \mathbb{E}\left[\xi^4_{n+1}|\mathcal{F}_n\right] \leq 16.
\end{equation}


Most of the asymptotic analysis will be conducted by the increasing process of martingale $(M_n)$. In particular, we define the predictable quadratic variation  $\langle M \rangle_n$, given for all $n \geq 1$, by
\begin{eqnarray}
\label{quadvar}
\langle M \rangle_n &= &\sum_{k=1}^n \mathbb{E}[\Delta M_k^2| \mathcal{F}_{k-1}] =\mathbb{E}\left[\xi^2_{1}|\mathcal{F}_0 \right]+ \sum_{k=1}^{n-1} a_{k+1}^2 \mathbb{E}\left[\xi^2_{k+1}|\mathcal{F}_k\right] \notag\\
&=& 1-2\omega(2\beta -1)+\omega^2 + \gamma\sum_{k=1}^{n-1} a_{k+1}^2 \frac{Z_k}{k}+(\tau-\omega^2)v_n\notag\\
&-&2\omega \alpha \sum_{k=1}^{n-1} a_{k+1}^2\frac{S_k}{k}-\alpha^2 \sum_{k=1}^{n-1} a_{k+1}^2 \frac{S_k^2}{k^2} \label{procrec}
\end{eqnarray}
where
\begin{equation} \label{vn}
v_n=\sum_{k=1}^n a_k^2.
\end{equation}
From well known results on the asymptotics of the gamma function, we may conclude that

\begin{enumerate}
\item If $\alpha<1/2$ then 
\begin{equation}
\label{vn1}
\frac{v_n}{n^{1-2\alpha}}\rightarrow \ell=\frac{\Gamma^2(\alpha+1)}{1-2\alpha}
\end{equation}

\item If $\alpha=1/2$ then
\begin{equation}
\label{vn2}
\frac{v_n}{\log n}\rightarrow \frac{\pi}{4}
\end{equation}

\item If $\alpha>1/2$ then, from \eqref{an} it is possible to deduce that $(v_n)$ converges into a finite value, more precisely
\begin{equation} \label{vn3}
v_n \rightarrow \sum_{k=0}^ \infty \left(\frac{\Gamma(\alpha+1)\Gamma(k+1)}{\Gamma(k+\alpha+1)} \right)^2 = \setlength\arraycolsep{1pt}
{}_3 F_2\left(\begin{matrix}1& &1& &1\\&\alpha+1&
&\alpha+1&\end{matrix};1\right),
\end{equation} 
where the above limit is the generalized hypergeometric function.
\end{enumerate}

These convergences will be the core for demonstrating the asymptotic results for the LRRW. On this direction, let us define the explosion coefficient of the martingale $M_n$ for $n\geq 1$ by
\begin{equation} \label{fnazero}
f_n=\frac{a_n^2}{v_n}
\end{equation}


Moreover, we will use adaptions of Theorems A.1 and A.2 from \cite{gue}, their original versions may be found in \cite{Cha} and \cite{cvgm} respectively.

\begin{lem}
\label{Lema1}
Let $\Delta M_n = M_n - M_{n-1}$, assume for all $\varepsilon > 0$

\begin{equation}
\label{cond1}
    \displaystyle\sum_{n=1}^{\infty} \frac{1}{v_n} \dE\left[\vert\Delta M_n\vert^2 \mathbb{I}_{\{\vert\Delta M_n\vert\geq\varepsilon\sqrt{v_n}\}}\vert\mathcal{F}_{n-1}\right] < \infty \ \  a.s
\end{equation}
and for some $a>0$,
\begin{equation}
\label{cond2}
    \displaystyle\sum_{n=1}^{\infty} \frac{1}{v_n^a} \dE\left[\vert\Delta M_n\vert^{2a} \mathbb{I}_{\{\vert\Delta M_n\vert\leq\sqrt{v_n}\}}\vert\mathcal{F}_{n-1}\right] < \infty \ \  a.s
\end{equation}
Then, $(M_n)$ satisfies that
\begin{equation}
  \displaystyle\frac{1}{\log v_{n}}\sum_{k=1}^n \left(\frac{v_{k}-v_{k-1}}{v_{k}}\right) \delta_{M_k/\sqrt{v_{k-1}}} \Rightarrow G  \ \  \text{ a.s}
\end{equation}
where $G$ stands for the $N(0,\sigma^2)$ distribution.
\end{lem}

\begin{lem}
\label{Lema2}
Assume that $(\xi_n)$, satisfies for some integer $r\ge 1$ and for some real number $h > 2r$ that
\begin{equation}
\label{cond3}
    \displaystyle\sup_{n\ge 0}\dE\left[\vert\xi_{n+1}\vert^h\vert\mathcal{F}_{n}\right] < \infty \ \  a.s
\end{equation}
If in addition the explosion coefficient $f_n$ tends a.s. to zero, then
\begin{equation}
    \displaystyle\lim_{n\to\infty} \frac{1}{\log v_n}\sum_{k=1}^n f_k \left( \frac{M_k^2}{v_{k-1}} \right)^r = \frac{\sigma^{2r} (2r)!}{2^r r!} \ \  a.s
\end{equation}
\end{lem}

%

Finally, we state the following lemma without its proof. However it may be obtained from Theorem 1 and Corollaries 1 and 2 from \cite{heyde1977central}.

\begin{lem}
\label{Heydelem}
Suppose that $\{M_n\}$ is a square-integrable martingale with mean $0$. Let $\Delta M_k = M_k - M_{k-1}$, for $k=1, 2, \ldots$, where $M_0=0$. If in addition
\begin{equation*}
    \displaystyle\sum_{k=1}^{\infty} \dE[(\Delta M_k)^2] < +\infty
\end{equation*}
holds, then we have the following: let $r_n^2= \displaystyle\sum_{k=n}^{\infty} \dE[(\Delta M_k)^2] $
\begin{itemize}
    \item [(i)] The limit $M:= \sum_{k=1}^{\infty} \Delta M_k $ exists almost surely and $M_n \overset{L^2}{\rightarrow} M$
    \item [(ii)] Assume that
    \begin{itemize}
        \item [a)] $\displaystyle\frac{1}{r_n^2}\sum_{k=n}^{\infty}\dE\left[(\Delta M_k)^2\vert\mathcal{F}_{k-1}\right]\rightarrow 1$ as $n \to \infty$ in probability, and
        \item [b)] $\displaystyle\frac{1}{r_n^2}\displaystyle\sum_{k=n}^{\infty}\dE\left[(\Delta M_{k+1})^2: \vert\Delta M_{k+1}\vert\geq\varepsilon r_k\right] \rightarrow 0$ for any $\varepsilon > 0$. \\
        Then, we have
        \begin{equation*}
            \displaystyle\frac{M-M_n}{r_{n+1}}=\frac{\sum_{k=n+1}^{\infty}\Delta M_k}{r_{n+1}} \overset{d}{\rightarrow}N(0,1).
        \end{equation*}
    \end{itemize}
     \item [(iii)] Assume that the following four conditions hold 
       \begin{itemize}
        \item [a')] $\displaystyle\frac{1}{r_n^2}\sum_{k=n}^{\infty}\dE\left[(\Delta M_k)^2\vert\mathcal{F}_{k-1}\right]\rightarrow 1$ as $n \to \infty$ a.s,
        \item [c)] $\displaystyle\sum_{k=1}^{\infty}\frac{1}{r_k}\dE[ |\Delta M_k| : | \Delta M_k|  > \varepsilon r_k] < +\infty$ for any $\varepsilon > 0$,
        
        \item [d)] $\displaystyle\sum_{k=1}^{\infty}\frac{1}{r_k^4}\dE[ (\Delta M_k)^4 : |\Delta M_k|  \le \delta r_k] < +\infty$ for some $\delta > 0$, and
            \item [e)] $\displaystyle\sum_{k=1}^{\infty}\frac{1}{r^2_k} \left( (\Delta M_k)^2 - \dE[ (\Delta M_k)^2 \vert\mathcal{F}_{k-1}]\right) < +\infty$~ a.s. 
    \end{itemize}
    Then, $\displaystyle\limsup_{n\to\infty}\pm \frac{M-M_n}{\sqrt{2\cdot r^2_{n+1} \log |\log r^2_{n+1} |}} = 1$ a.s.
\end{itemize}
\end{lem}



\subsection{The diffusive regime}
\subsubsection{Proof of theorem \ref{CVG1}}

From \eqref{procrec} and \eqref{vn1} we note that
\begin{equation}
\langle M \rangle_n = O(v_n)=O(n^{1-2\alpha}),
\end{equation}
which together with Theorem $1.3.24$ of \cite{Duflo} imply $|M_n|=O\left(\sqrt{n^{1-2\alpha} \log n }\right)$. By recalling that,  $M_n=a_n S_n-\omega A_n$ and that $a_n \sim \frac{\Gamma(1+\alpha)}{n^\alpha}$, we notice that
\begin{equation} \label{dif1}
\left|\frac{S_n}{n}-\frac{\omega A_n}{n a_n}\right|=O\left(\sqrt{\frac{\log n}{n} }\right).
\end{equation}
However, Lemma $B.1$ of \cite{ERWBercu} implies that
\begin{equation} 
\frac{A_n}{n a_n}=\frac{1}{\alpha-1}\left(\frac{\Gamma(n+\alpha)}{\Gamma(n+1)\Gamma(\alpha)}-1\right), \label{fromlemB1}
\end{equation}
hence

\begin{equation}
\left|\frac{A_n}{n a_n}-\frac{1}{1-\alpha}\right| \sim  \frac{1}{(1-\alpha)\Gamma(\alpha)}\frac{1}{n^{1-\alpha}}, \label{eq233}
\end{equation}
which together with \eqref{dif1} implies that, almost surely $\frac{S_n}{n}\rightarrow \frac{\omega}{1-\alpha}.$ More accurately, from \eqref{dif1} and \eqref{eq233} we obtain that:
\begin{eqnarray}
\Big(\frac{S_n}{n}-\frac{\omega}{1-\alpha}\Big)^2\!\!&=&\!\Big(\frac{S_n}{n}-\frac{\omega A_n}{n a_n}\Big)^2\!\!+\!\Big(\frac{\omega A_n}{n a_n}-\frac{\omega}{1-\alpha}\Big)^2\\
&+&2\Big(\frac{S_n}{n}-\frac{\omega A_n}{n a_n}\Big)\!\Big(\frac{\omega A_n}{n a_n}-\frac{\omega}{1-\alpha}\Big) \notag \\
&=&O\left( \frac{\log n}{n} \right).\notag
\end{eqnarray}

That concludes the proof of Theorem \ref{CVG1}.


\subsubsection{Proof of Theorem \ref{thmmoments}}
Condition \eqref{cond3} follows; for all $m \in \mathbb{N}$, since $\dE\left(X_{n+1}^{m} | \mathcal{F}_n\right) \leq 1,$ which conduces us to
\begin{eqnarray*}
\dE[(\xi_{n+1})^m\vert\mathcal{F}_{n}] & = & \dE\left[\left(X_{n+1}-\dE(X_{n+1}\vert \mathcal{F}_{n})\right)^m\vert\mathcal{F}_{n}\right] \\
                                      & = & \dE\left[\left.\sum_{i=0}^{m}{m \choose i}\left(X_{n+1}\right)^{m-i}\left(\dE(X_{n+1}\vert \mathcal{F}_{n})\right)^i(-1)^i \ \right| \ \mathcal{F}_{n}\right] \\
                                      & \leq & (m+1) \max_{i=0,\cdots,h}{m \choose i} \,
\end{eqnarray*}
which implies that $\sup_{n \geq 0}\dE\left[|\xi_{n+1}|^h \vert \mathcal{F}_n\right] < + \infty$. Hence, it is only needed to check that $f_n$ goes to zero as $n\to\infty$ by virtue of Lemma \ref{Lema2}. Nevertheless, it follows from definition \eqref{fnazero} and convergences \eqref{an} and \eqref{vn1}, which imply that $f_n\sim\frac{1-2\alpha}{n}$.

 Finally, by noticing that  $\frac{M_k}{\sqrt{v_{k-1}}}  \sim \sqrt{(1-2\alpha)k} \left( \frac{S_k}{k} -\frac{\omega}{1-\alpha} \right) \ \text{a.s.}$ the Theorem holds.

\subsubsection{Proof of Theorem \ref{CVG3}}

In order to demonstrate Theorem \ref{CVG3} we analyse the behaviour of $Z_n$. Then, by \eqref{xsqq} we observe that $\mathbb{E}\left[Z_{n+1}|\mathcal{F}_n\right]=\lambda_n Z_n+\tau \ \text{a.s.,}$ where $\lambda_n=1+\frac{\gamma}{n}$. Which leads us to consider the discrete time martingale $(N_n)$ given, for $n\geq1$ by $N_n=b_n Z_n-\tau B_n,$ where $b_n=\frac{\Gamma(n)\Gamma(\gamma+1)}{\Gamma(n+\gamma)}$ and $B_n=\displaystyle\sum_{k=1}^n b_k$,
for proving (from similar arguments of those presented in proof of Theorem \ref{CVG1}) that $\frac{Z_n}{n}\rightarrow \frac{\tau}{1-\gamma} \ \text{a.s.}$, which jointly with \eqref{segunda} and Theorem \ref{CVG1} lets us to see that
\begin{equation}\label{sigma222}
\lim_{n\rightarrow \infty} \mathbb{E}\left[\xi^2_{n+1}|\mathcal{F}_n\right] = \sigma^2 \hspace{.5cm}\text{a.s.,}
\end{equation}
where $\sigma^2$ as in \eqref{notation}. Moreover, note that \eqref{vn1} implies 
$\sum_{k=1}^\infty \frac{a_k^4}{v_n^2}=\frac{\left(\pi (1-2\alpha\right))^2}{6}.$
Finally, the law of iterated logarithm due to Stout \cite{stout} (or the one in Lemma C.2 of \cite{bblil}) leads us to affirm that almost surely

\begin{equation} \label{lil11}
\limsup_{n\rightarrow \infty} \frac{M_n}{\sqrt{2v_n\log \log v_n}}=-\liminf_{n\rightarrow \infty} \frac{M_n}{\sqrt{2v_n\log \log v_n}}=\sigma.
\end{equation}

From definition of $(M_n)$ and convergence \eqref{vn1}, we obtain that
\begin{eqnarray} 
&\limsup_{n\rightarrow \infty}& \left( \frac{n}{2\log \log n} \right)^{1/2} \left(\frac{S_n}{n}-\frac{\omega A_n}{n a_n}\right) \notag\\
  &=&-\liminf_{n\rightarrow \infty} \left( \frac{n}{2\log \log n} \right)^{1/2} \left(\frac{S_n}{n}-\frac{\omega A_n}{n a_n}\right) \notag\\
&=&\frac{\sigma}{\sqrt{1-2\alpha}}. \notag
\end{eqnarray}
However, from \eqref{fromlemB1} and definition of $a_n$, we have that
\begin{eqnarray} 
&\limsup_{n\rightarrow \infty}& \left( \frac{n}{2\log \log n} \right)^{1/2} \left(\frac{S_n}{n}- \frac{\omega}{1-\alpha}-\frac{\alpha \omega}{(\alpha-1)na_n} \right) \notag\\
  &=&-\liminf_{n\rightarrow \infty} \left( \frac{n}{2\log \log n} \right)^{1/2} \left(\frac{S_n}{n}-\frac{\omega}{1-\alpha}-\frac{\alpha \omega}{(\alpha-1)na_n}\right) \notag\\
&=&\frac{\sigma}{\sqrt{1-2\alpha}},\notag
\end{eqnarray}
But, given that $\alpha<1/2$ we find that, as $n \to \infty $ $\left( \frac{n}{2\log \log n} \right)^{1/2}\left(\frac{1}{n^{1-\alpha}}\right)\rightarrow 0$, which implies the conclusion of Theorem \ref{CVG3}.

\subsubsection{Proof of Theorem \ref{thmASCLT}}

The proof is essentially based on Lemma \ref{Lema1}. Hence, note that via \eqref{martingale2}, we have that

\begin{equation*}
     \begin{array}{l}
    \displaystyle\displaystyle\sum_{k=1}^{\infty}\frac{1}{v_k}\dE\left[\vert\Delta M_{k}\vert^2 \mathbb{I}_{\vert\Delta M_{k}\vert\geq\varepsilon\sqrt{v_k}}\vert\mathcal{F}_{k-1}\right]       \leq \frac{1}{\varepsilon^2}\displaystyle\sum_{k=1}^{\infty}\frac{1}{v_k^2}\dE\left[\vert\Delta M_{k}\vert^4 \vert\mathcal{F}_{k-1}\right] \\[0.5cm]
           \leq \displaystyle\sup_{k\geq 1}\dE\left[\xi^4_{k} \vert\mathcal{F}_{k-1}\right]\frac{1}{\varepsilon^2}\displaystyle\sum_{k=1}^{\infty} \frac{a_k^4}{v_k^2 } \leq\frac{16}{\varepsilon^2}\displaystyle\sum_{k=1}^{\infty} \frac{ a_k^4}{v_k^2} \sim \frac{16}{\varepsilon^2}\displaystyle\sum_{k=1}^{\infty} \frac{(1-2\alpha)^2}{k^2} < \infty
     \end{array}
 \end{equation*}
Where, last step is due to \eqref{an} and \eqref{vn1}. Therefore \eqref{cond1} holds. To prove the validity of \eqref{cond2} we follow analogous steps with $a=2$.
Then
\begin{equation}
  \displaystyle\frac{1}{\log v_{n}}\sum_{k=1}^n \left(\frac{v_{k}-v_{k-1}}{v_{k}}\right) \delta_{M_k/\sqrt{v_{k-1}}} \Rightarrow G^*  \ \  \text{ a.s}
\end{equation}
By recalling that, $f_k \sim \frac{1-2\alpha}{k}$, $\log v_n \sim (1-2\alpha) \log n$ and $\frac{M_k}{\sqrt{v_{k-1}}} \sim \sqrt{\frac{1-2\alpha}{k}} \left(S_k - k\frac{\omega }{1-\alpha}\right)$, we conclude that
\begin{equation}
  \displaystyle\frac{1}{\log n}\sum_{k=1}^n \frac{1}{k} \delta_{\sqrt{k}\left(\frac{S_k}{k} - \frac{\omega }{1-\alpha}\right)} \Rightarrow G  \ \  \text{ a.s}
\end{equation}
where $G \sim N(0, \sigma^2 /(1-2 \alpha))$ is the re-scaled version of $G^* \sim N(0,\sigma^2)$. 
\\

\subsubsection{Proof of Theorem \ref{T-FCLT-DR}}
Note that, using \eqref{quadvar} and Toeplitz lemma \cite{Duflo}, we have

\begin{eqnarray*}
\lim_{n\rightarrow \infty} \frac{1}{n^{1-2\alpha}}\langle M \rangle_n &=&\displaystyle\frac{\Gamma(\alpha+1)^2}{1-2\alpha}\left( \frac{\gamma\tau}{1-\gamma} + (\tau-\omega^2) -  \frac{2\omega^2\alpha}{1-\alpha} - \left(\frac{\omega\alpha}{1-\alpha}\right)^2 \right)\\
& =& \sigma^2\ell \hspace{1cm}\text{a.s.}
\label{CVGIP-DR}
\end{eqnarray*}
where $\ell$ is given in \eqref{vn1}. Then, we apply the functional central limit theorem for martingales given in Theorem 2.5 of \cite{DR}. That is, consider the martingale difference array $D_{n,k} = \frac{1}{\sqrt{n^{1-2\alpha}}}(\Delta M_k)$, which satisfies
\begin{equation}
\label{convfunctional}
\lim_{n\rightarrow \infty} \frac{1}{n^{1-2\alpha}}\langle M \rangle_{\lfloor nt \rfloor} = \sigma^2\ell t^{1-2\alpha} \hspace{1cm}\text{a.s.}
\end{equation}
In addition, we need to prove the Lindeberg's condition. We obtain from \eqref{cuarta} that for any $\varepsilon > 0$
\begin{eqnarray*}
\displaystyle\frac{1}{n^{1-2\alpha}}\sum_{k=1}^n\dE [\Delta M_k ^2 \mathbb{I}_{\{|\Delta M_k| > \varepsilon \sqrt{n^{1-2\alpha}}\}} \vert \mathcal{F}_{k-1}] &\leq&\displaystyle\frac{1}{n^{2(1-2\alpha)}\varepsilon^2}\sum_{k=1}^n\dE [\Delta M_k ^4 \vert \mathcal{F}_{k-1}]\\
\leq \displaystyle\frac{1}{n^{2(1-2\alpha)}\varepsilon^2}\sum_{k=1}^n a_k^4\dE [\xi_k ^4 \vert \mathcal{F}_{k-1}]
 &\leq& \displaystyle\frac{16}{n^{2(1-2\alpha)}\varepsilon^2}\sum_{k=1}^n a_k^4,
\end{eqnarray*}
Then, thanks to \eqref{an} and \eqref{vn1}, we have that as $n \rightarrow \infty$, $\frac{n^2 a_n^4}{v_n^2}\rightarrow (1-2\alpha)^2$, which implies that $\frac{1}{n^{1-4\alpha}}\sum_{k=1}^n a_k^4$ converges to $\frac{(1-2\alpha)^2\ell^2}{1-4\alpha}$. Therefore,
\begin{equation*}
\displaystyle\frac{1}{n^{1-2\alpha}}\sum_{k=1}^n\dE [\Delta M_k ^2 \mathbb{I}_{\{|\Delta M_k| > \varepsilon \sqrt{n^{1-2\alpha}}\}} \vert \mathcal{F}_{k-1}] \rightarrow 0 \text{ as } n \to \infty \text{ in probability},
\end{equation*}
which allows us to conclude that for all $t\ge 0$ and for any $\varepsilon > 0$, 
\begin{equation}
\label{Lindfunctional}
\displaystyle\frac{1}{n^{1-2\alpha}}\sum_{k=1}^{\lfloor nt \rfloor}\dE [\Delta M_k ^2 \mathbb{I}_{\{|\Delta M_k| > \varepsilon \sqrt{n^{1-2\alpha}}\}} \vert \mathcal{F}_{k-1}] \rightarrow 0,
\end{equation}
as $n \to \infty$ in probability. By noticing that $\lim_{n\to\infty}\frac{{\lfloor nt \rfloor}a_{\lfloor nt \rfloor}}{n^{1-2\alpha}}=t^{1-\alpha} \Gamma(\alpha+1)$ and that \eqref{an} and \eqref{fromlemB1} imply that
\begin{equation}
\label{Mntfunctional}
 \frac{M_{\lfloor nt \rfloor}}{\sqrt{n^{1-2\alpha}}}  = \frac{{\lfloor nt \rfloor}a_{\lfloor nt \rfloor}}{\sqrt{n^{1-2\alpha}}} \left( \frac{S_{\lfloor nt \rfloor}}{{\lfloor nt \rfloor}} -\frac{\omega}{1-\alpha}\right) + \frac{\omega\alpha}{(1-\alpha)\sqrt{n^{1-2\alpha}}} \hspace{1cm}\text{a.s.,}
\end{equation}
we conclude via Theorem 2.5 of \cite{DR} that $\left( \sqrt{n} \left( \frac{S_{\lfloor nt \rfloor}}{{\lfloor nt \rfloor}} -\frac{\omega}{1-\alpha}\right), t \geq 0\right) \Longrightarrow \big( W_t, t \geq 0 \big),$ where $W_t = B_t /(t^{1-\alpha} \Gamma(\alpha+1))$, which completes the proof of the theorem.

\subsection{The critical regime}

\subsubsection{Proof of Theorem \ref{CVG11}}

We will proceed in a similar fashion as in Theorem \ref{CVG1}. First of all, from the last part of Theorem 1.3.24 of \cite{Duflo} we have that $\frac{M_n^2}{\log n}=O(\log \log n)\hspace{.2cm}\text{a.s.}$
In addition, \eqref{an} provides us that
\begin{equation}\label{eq555}
na_n^2\rightarrow \frac{\pi}{4}.
\end{equation}
Moreover, from the definition of $(M_n)$ we may see that $\left|a_n S_n-\omega A_n\right|=O(\sqrt{\log n \log \log n}) \ \text{a.s}$. Hence, from \eqref{eq555}, we notice that
\begin{equation} \label{llncrit}
\left|\frac{S_n}{n}-\frac{\omega A_n}{na_n}\right|=O\left(\sqrt{\frac{\log n \log \log n}{n}}\right)\hspace{.2cm}\text{a.s.}
\end{equation}
Additionally, we observe from \eqref{an} that $A_n\sim \sqrt{n\pi},$ that, together with \eqref{eq555} let us notice that $\frac{A_n}{na_n}\rightarrow 2.$ Furthermore, from analogue lines than those which imply \eqref{eq233} we obtain that
\begin{equation}\label{eq4433}
\left|\frac{A_n}{n a_n}-2\right|\sim \frac{2}{\sqrt{n\pi}}. 
\end{equation}
Which, together with \eqref{llncrit} and the triangle inequality leads us to \eqref{as001}, and the proof of Theorem \ref{CVG11} is achieved. 

\subsubsection{Proof of Theorem \ref{thm4}}

This theorem may be proven in a similar fashion than Theorem \ref{thmmoments} and Theorem 2.3 of \cite{BV2021}. For this, we notice that condition \eqref{cond3} holds in the same manner than in the diffusive regime. Besides, we note that Theorem \ref{CVG11} together with \eqref{segunda} implies the following almost sure convergence

\begin{equation}
\lim_{n\rightarrow \infty} \mathbb{E}\left[\xi^2_{n+1}|\mathcal{F}_n\right] = \frac{\tau}{1-\gamma}-4\omega^2 .
\end{equation}

In addition, due to \eqref{an} and \eqref{vn2} we obtain that 
\begin{equation}\label{ll2}
\frac{a_k^4}{v_k^2}\sim \left(\frac{1}{n \log n}\right)^2,
\end{equation} %
which implies that $f_n$ converges to zero as $n \rightarrow \infty$. Hence, we may conclude  \eqref{moments2}  from the definition of $(M_n)$, and Theorem $3$ of \cite{cvgm}.


\subsubsection{Proof of Theorem \ref{CVG31}}

This theorem follows in a similar way than Theorem \ref{CVG3}, by noticing that equation \eqref{ll2} implies that

\begin{equation} \label{ll3}
\sum_{k=1}^\infty \frac{a_k^4}{v_n^2}<\infty.
\end{equation}



\subsubsection{Proof of Theorem \ref{thmasclt}}

Note that the conditions of Lemma \ref{Lema1} follows from \eqref{ll3}. In addition, it may be found from the definition of $M_n$, \eqref{vn2}, \eqref{eq4433} and \eqref{ll2} that $\frac{M_n}{\sqrt{v_{n-1}}} \sim \sqrt{\frac{n}{\log n}}\left[ \frac{S_n}{n}-2\omega \right],$ which leads to $\displaystyle\frac{1}{\log \log n}\sum_{k=1}^n \frac{1}{k \log k} \delta_{\sqrt{\frac{k}{\log k}}\left(\frac{S_k}{k} - 2\omega \right)} \Rightarrow G  \ \text{ a.s.,}$ where $G$ stands for the $N(0,\sigma^2)$ distribution.



\subsubsection{Proof of Theorem \ref{T-FCLT-CR}}

This proof is essentially the same that of Theorem \ref{T-FCLT-DR}, by noticing that Lindeberg condition holds from \eqref{ll3}. Then, from the functional central limit theorem for martingales \cite{DR}, the definition of $(M_n)$, convergence \eqref{vn2}, and relation equation \eqref{ll2} we reach the Theorem.

\subsection{The superdiffusive regime}
\subsubsection{ Proof of Theorem \ref{T-ASP-SR}.}

From decomposition \eqref{martingale2}, \eqref{epsn} and \eqref{cota112}, we have; for $n\geq 0$, that

$$\sup_{n\geq 1} \mathbb{E}[M_{n}^2]\leq (\gamma + \tau) \cdot \setlength\arraycolsep{1pt}
{}_3 F_2\left(\begin{matrix}1& &1& &1\\&\alpha+1&
&\alpha+1&\end{matrix};1\right) <\infty,$$
since $v_n$ is a non decreasing sequence. That is to say, martingale $(M_n)$ is bounded in $L^2$. Thus, it converges in $L^2$ and almost surely to the random variable
\begin{equation}\label{limitM}
M=\sum_{k=1}^\infty a_k \xi_k.
\end{equation}
The almost sure convergence $M_n=a_n S_n-\omega A_n \rightarrow M$ and \eqref{fromlemB1} infer that $na_n\left( \frac{S_n}{n}-\frac{\omega}{1-\alpha}\right) \rightarrow M+\frac{\omega\alpha}{\alpha-1}\hspace{.2cm}\text{a.s.,}$ which, jointly with \eqref{an} implies that $n^{1-\alpha}\left( \frac{S_n}{n}-\frac{\omega}{1-\alpha}\right) \rightarrow L \hspace{.2cm}\text{a.s.,}$ where
\begin{equation}\label{defL}
L=\frac{1}{\Gamma(\alpha+1)}\left(M-\frac{\omega\alpha}{1-\alpha}\right).
\end{equation}
Moreover, given that $(M_n)$ converges to M in $L^2$ we have that \eqref{fromlemB1} and definition of the limit random variable $L$ guide us to
$$\lim_{n\rightarrow \infty} \mathbb{E}\left[\left(n^{1-\alpha}\left(\frac{S_n}{n}-\frac{\omega}{1-\alpha}\right) -L \right)^2 \right]=0,$$
We will find now the first two moments of the limiting random variable $L$ given in \eqref{defL}. For this, let us note that $\mathbb{E}[X_1]=\beta$. In addition, from \eqref{expos} we have; for $n=1,2,\ldots,$ that
\begin{equation}
\label{lEsn}
    \dE[S_n]=\frac{1}{ a_n}\left(\beta+\omega\cdot  \displaystyle\sum_{l=1}^{n-1}a_{l+1}\right) = \frac{1}{ a_n}\left(\beta+\omega\cdot (A_n-1)\right)
\end{equation}
we directly see that $ \mathbb{E}[M_n]=\beta-\omega,$ which implies that $\lim_{n\rightarrow \infty} \mathbb{E}[M_n]=\mathbb{E}[M]=\beta-\omega,$ that leads us to
$$\mathbb{E}(L)=\frac{1}{\Gamma(\alpha+1)}\left(\mathbb{E}(M)-\frac{\omega\alpha}{1-\alpha}\right)=\frac{\beta(1-\alpha)-\omega}{\Gamma(\alpha+1)(1-\alpha)},$$
that implies \eqref{meanL}. We will proceed now with the deduction of \eqref{varL}. For this, note that
\begin{eqnarray}
\mathbb{E}[M_n^2]&=&a_n^2\mathbb{E}[S_n^2]-2\omega a_n A_n\mathbb{E}[S_n]+\omega ^2A_n^2 \notag\\[0.1cm]
&=&a_n^2\mathbb{E}[S_n^2]-2\omega  A_n(\beta-\omega)-\omega ^2A_n^2. \label{mnsquaree}
\end{eqnarray}
Moreover $\mathbb{E}[S_{n+1}^2]=g_n\mathbb{E}[S_{n}^2]+h_n$, where, $g_n:=1+\frac{2\alpha}{n}$, and $h_n:=2\omega\mathbb{E}[S_n]+\frac{\gamma}{n}\dE(Z_n)+\tau$, for $n \geq 1$. Hence, it may be found recursively; for $n \geq 1$, that
\begin{equation} \label{snsquare}
\mathbb{E}[S_{n}^2]=\frac{\Gamma(n+2\alpha)}{\Gamma(n)\Gamma(2\alpha+1)}\left(p+q+\Gamma(2\alpha+1)\sum_{k=1}^{n-1}h_k\frac{\Gamma(k+1)}{\Gamma(k+1+2\alpha)}\right).
\end{equation}
However, from \eqref{fromlemB1} we may see that
\begin{eqnarray*}
h_k&=&2\omega \left[ \frac{\beta-\omega(1-A_k)}{a_k}\right]+\frac{\gamma}{k}\left[ \frac{p+q-\tau(1-B_k)}{b_k} \right]+\tau \\
&=&\frac{2\omega(\beta-\omega)}{a_k}+2\omega^2\left[ \frac{k}{1-\alpha}-\frac{\alpha}{(1-\alpha)a_k}\right]+\frac{\gamma(p+q-\tau)}{kb_k}\\
&+&\gamma \tau \left[\frac{1}{1-\gamma}-\frac{\gamma}{(1-\gamma)kb_k}+\frac{1}{\gamma}\right] = \frac{\tau}{1-\gamma}+\frac{2k\omega^2}{1-\alpha}-\frac{t_1}{a_k}+\frac{\gamma}{k b_k}t_2,
\end{eqnarray*}
where $t_1:=\frac{2\omega[\alpha \omega-(1-\alpha)\beta \theta]}{1-\alpha}$ and $t_2:=\gamma r$. In addition,
\begin{eqnarray*}
h_k\frac{\Gamma(k+1)}{\Gamma(k+1+2\alpha)} &=& \frac{\tau}{1-\gamma}\frac{\Gamma(k+1)}{\Gamma(k+1+2\alpha)}+\frac{2\omega^2}{1-\alpha}\frac{k\Gamma(k+1)}{\Gamma(k+1+2\alpha)}\\
&-&t_1\frac{\Gamma(k+1)}{a_k\Gamma(k+1+2\alpha)}+\gamma t_2\frac{\Gamma(k+1)}{k b_k\Gamma(k+1+2\alpha)}
\end{eqnarray*}

Through a very exhaustive employ of Lemma $B.1$ of \cite{ERWBercu} we observe that
\begin{eqnarray*}
\sum_{k=1}^{n-1} h_k \frac{\Gamma(k+1)}{\Gamma(k+1+2\alpha)} &=&\frac{\tau}{1-\gamma}\left[\frac{1}{(2\alpha-1)\Gamma(1+2\alpha)}-\frac{n\Gamma(n)}{(2\alpha-1)\Gamma(n+2\alpha)}\right]\\
&+&\frac{\omega^2}{1-\alpha}\left[ \frac{1}{(\alpha-1)(2\alpha-1)\Gamma(2\alpha)}-\frac{n\Gamma(n)\left(n(2\alpha-1)+1\right)}{(\alpha-1)(2\alpha-1)\Gamma(n+2\alpha)}\right] \\
&-&t_1 \left[\frac{1}{\alpha(\alpha-1)\Gamma(2\alpha)}-\frac{(n+1)\Gamma(n+\alpha)}{\alpha(\alpha-1)\Gamma(\alpha)\Gamma(n+2\alpha)}\right]\\
&+&\frac{\gamma t_2}{\Gamma(\gamma+1)}\left[ \frac{\Gamma(\gamma+1)}{(2\alpha-\gamma)\Gamma(2\alpha+1)}-\frac{\Gamma(n+\gamma)}{(2\alpha-\gamma)\Gamma(2\alpha+n)} \right].
\end{eqnarray*}

From this, after several arithmetical calculations, thanks to \eqref{snsquare} we find that
\begin{eqnarray*}
\mathbb{E}[S_n^2]&=&\frac{\Gamma(n+2\alpha)}{\Gamma(n)\Gamma(2\alpha+1)}\nabla-\frac{n}{2\alpha-1}\left[\frac{\tau}{1-\gamma}-\frac{\omega^2\left(n(2\alpha-1)+1\right)}{(\alpha-1)^2} \right]\\
&+&\frac{(n+1)t_1}{a_n(\alpha-1)}-\frac{\gamma t_2}{(2\alpha-\gamma)b_n},
\end{eqnarray*}
where $\nabla$ is given by \eqref{varL}. Hence, from \eqref{mnsquaree}, we imply that
\begin{eqnarray*}
\mathbb{E}[M_n^2]&=&\frac{a_n^2 \Gamma(n+2\alpha)\nabla}{\Gamma(n)\Gamma(2\alpha+1)}-\frac{\tau n a_n^2}{(2\alpha-1)(1-\gamma)}+\frac{n \omega^2 a_n^2}{(2\alpha-1)(\alpha-1)^2}\\
&-&\frac{\gamma t_2 a_n^2}{(2\alpha-\gamma)b_n}-2\omega[\beta-\omega]\frac{n a_n-\alpha}{1-\alpha}+\frac{2\alpha\omega^2 n a_n}{(1-\alpha)^2}\\
&-&\frac{\alpha^2 \omega^2}{(1-\alpha)^2}+\frac{t_1(n+1)a_n}{\alpha-1},
\end{eqnarray*}
however, it may be obtained that $\frac{t_1}{\alpha-1}+\frac{2\alpha \omega^2}{(1-\alpha)^2}-\frac{2\omega[\beta-\omega]}{1-\alpha}=0$, which leads us to
\begin{eqnarray*}
\mathbb{E}[M_n^2]&=&\frac{a_n^2 \Gamma(n+2\alpha)\nabla}{\Gamma(n)\Gamma(2\alpha+1)}-\frac{\tau n a_n^2}{(2\alpha-1)(1-\gamma)}+\frac{n \omega^2 a_n^2}{(2\alpha-1)(\alpha-1)^2}\\
&-&\frac{\gamma t_2 a_n^2}{(2\alpha-\gamma)b_n}+2\omega[\beta-\omega]\frac{\alpha}{1-\alpha} -\frac{\alpha^2 \omega^2}{(1-\alpha)^2}+\frac{t_1 a_n}{\alpha-1},
\end{eqnarray*}
that is such that $\mathbb{E}[M^2]=\lim_{n\rightarrow \infty} \mathbb{E}[M_n^2]=\frac{\Gamma^2(\alpha+1) \nabla}{\Gamma(2\alpha+1)}+\frac{\alpha^2 \omega(2-\beta(\theta+1))}{(1-\alpha)^2}$. So, from \eqref{defL} and the first two moments of $M$ we conclude that
\begin{eqnarray}
\mathbb{E}[L^2]&=&\frac{1}{\Gamma^2(\alpha +1)}\left[\mathbb{E}[M^2]-\frac{2\omega \alpha}{1-\alpha}\mathbb{E}[M]+\frac{\omega^2 \alpha^2}{(1-\alpha)^2}\right]\notag \\
&=&\frac{\nabla}{\Gamma(2\alpha+1)}+2\omega\left(\frac{1}{(1-\alpha)\Gamma(\alpha)}\right)^2
\end{eqnarray}

\subsubsection{Proof of Theorem \ref{Gaussfluct}}

Note that, \eqref{sigma222} and the bounded convergence theorem imply that $\displaystyle\sum_{k=1}^{\infty}\dE[(\Delta M_k)^2]\sim \sigma^2\Gamma(\alpha+1)^2\displaystyle\sum_{k=1}^{\infty}\frac{1}{k^{2\alpha}}$. Then, since $\alpha> 1/2$, we have that $\displaystyle\sum_{k=1}^{\infty}\dE[(\Delta M_k)^2] < \infty$.

Now, given \eqref{sigma222}, we have that:

\begin{equation*}
\begin{array}{ll}
     \displaystyle\sum_{k=n}^{\infty}\dE[(\Delta M_k)^2\vert\mathcal{F}_{k-1}]&\sim\displaystyle\sum_{k=n}^{\infty} \sigma^2 a_k^2
    \sim \sigma^2\Gamma(\alpha+1)^2\displaystyle\sum_{k=n}^{\infty}\frac{1}{k^{2\alpha}}\\[0.4cm]
    &\sim\displaystyle\frac{\sigma^2\Gamma(\alpha+1)^2}{(2\alpha-1)n^{2\alpha-1}}
    \sim \frac{\sigma^2 }{(2\alpha-1)} n a_n^2 ~a.s
\end{array}
\end{equation*}
By using the bounded convergence theorem, we obtain that:
\begin{equation}
\label{r2conv}
    r_n^2 := \sum_{k=n}^{\infty}\dE[(\Delta M_k)^2]\sim \frac{\sigma^2 }{(2\alpha-1)}n a_n^2~a.s
\end{equation}
Then, conditions a) and a') in Lemma \ref{Heydelem} are satisfied. Hence, have that
\begin{equation*}
     r_n^4\sim \frac{\Gamma(\alpha+1)^4\sigma^4 }{(2\alpha-1)^2n^{4\alpha-2}} \text{ and } (\Delta M_n)^4\leq 16 a_n^4\sim \frac{16\Gamma(\alpha+1)^4}{n^{4\alpha}}.
\end{equation*}
In this sense, observe that
\begin{equation*}
\begin{array}{c}
         \displaystyle\frac{1}{r_n^2}\sum_{k=n}^{\infty}\dE\left[(\Delta M_{k+1})^2: \vert\Delta M_{k+1}\vert\geq\varepsilon r_k\right]  \leq \displaystyle\frac{1}{\varepsilon^2r_n^4}\displaystyle\sum_{k=n}^{\infty}\dE\left[\vert\Delta M_{k+1}\vert^4\right] \\[0,4cm]
 \leq \displaystyle\frac{16}{\varepsilon^2r_n^4}\sum_{k=n}^{\infty} a_k^4
         \leq\frac{C_1}{r_n^4}\displaystyle\sum_{k=n}^{\infty} \frac{1}{k^{4\alpha}}\leq \displaystyle C_2 n^{4\alpha-2}n^{1-4\alpha}
\end{array}
 \end{equation*}
 which implies condition b) of Lemma \ref{Heydelem}. Then, by noticing that $M_n-M = a_n \left(S_n - n\frac{\omega}{1-\alpha} -n^{\alpha}L\right)$, and using \eqref{r2conv}, the convergence \eqref{T_cnv_w} holds.

 Additionally, for $\varepsilon>0$ we have that:
\begin{equation*}
         \frac{1}{r_k}\dE\left[\vert\Delta M_{k+1}\vert: \vert\Delta M_{k+1}\vert\geq\varepsilon r_k\right] \leq\frac{1}{r_k} \frac{1}{\varepsilon^3r_k^3}\dE\left[\vert\Delta M_{k+1}\vert^4\right]
         \leq C_3 k^{4\alpha-2}k^{4\alpha}
         =\frac{C_3}{k^2}
 \end{equation*}
which implies that condition c) of the same Lemma is satisfied.
In addition, given that $\sum_{k=1}^{\infty} \frac{1}{r_k^4}\dE[(\Delta M_k)^4]<\infty$, we obtain the validation of condition d).

Finally, let us denote, $d_k:=\frac{1}{r_k^2} \left( (\Delta M_k)^2 - \dE[ (\Delta M_k)^2 \vert\mathcal{F}_{k-1}]\right)$, the martingale difference, and observe that
\begin{equation*}
\begin{array}{ll}
         \displaystyle\sum_{k=1}^{\infty}\dE\left[d_k^2 \vert\mathcal{F}_{k-1}\right]  = \displaystyle\sum_{k=1}^{\infty}\frac{1}{r_k^4}\left(\dE[ (\Delta M_k)^4 \vert\mathcal{F}_{k-1}]-\dE^2[ (\Delta M_k)^2 \vert\mathcal{F}_{k-1}]\right)\\[0,4cm]
        \leq \displaystyle\sum_{k=1}^{\infty}\frac{1}{r_k^4}\dE[ (\Delta M_k)^4 \vert\mathcal{F}_{k-1}]
         \leq C_4\displaystyle\sum_{k=1}^{\infty} \frac{a_k^4}{r_k^4}\leq C_5\displaystyle\sum_{k=1}^{\infty} \frac{1}{k^2} < + \infty
\end{array}
 \end{equation*}
 Then, as a consequence of Theorem 2.15 from \cite{Hall1980}, condition e) of Lemma \ref{Heydelem} is satisfied, and therefore \eqref{LIL3} holds.

\backmatter


\section*{Statements and Declarations}

\begin{itemize}
\item Funding: MGN was partially supported by Fondecyt Iniciaci\'on 11200500. RL is partially supported by CNPQ 435470/2018-3 and FAPEMIG APQ-01341-21 and RED-00133-21 projects. VHVG is partially supported by Proyectos VIEP-BUAP-00154.

\item Competing interests: Authors have no relevant financial or non-financial interests to disclose, nor any competing interests to declare that are relevant to the content of this article.

\end{itemize}






\bibliography{sn-bibliography}

\end{document}